\documentclass[12pt,oneside,reqno]{amsart}

\usepackage[letterpaper,textwidth=6.9in,textheight=9.6in,centering]{geometry}
\usepackage{amsmath, amsfonts, amssymb}
\usepackage{mathrsfs}
\usepackage{latexsym}
\usepackage{bbm}
\usepackage{mathtools}
\usepackage{exscale}
\usepackage{cases}

\usepackage{enumitem}

\usepackage[pdfstartview=FitH,
            CJKbookmarks=true,
            bookmarksnumbered=true,
            bookmarksopen=true,
            colorlinks=true,
            linkcolor=blue,
            anchorcolor=blue,
            citecolor=red,
            urlcolor=blue
            ]{hyperref}

\numberwithin{equation}{section}
\newtheorem{theorem}{Theorem}[section]

\newtheorem{lemma}[theorem]{Lemma}
\newtheorem{corollary}[theorem]{Corollary}

\newtheorem{proposition}[theorem]{Proposition}

\newtheorem{definition and theorem}[theorem]{Definition and Theorem}

\def\bl{\begin{lemma}}
\def\el{\end{lemma}}
\def\bc{\begin{corollary}}
\def\ec{\end{corollary}}
\def\bt{\begin{theorem}}
\def\et{\end{theorem}}

\def\bp{\begin{proposition}}
\def\ep{\end{proposition}}
\def\be{\begin{equation}}
\def\ee{\end{equation}}
\def\baa{\begin{align*}}
\def\eaa{\end{align*}}

\theoremstyle{definition}
 
\theoremstyle{remark}
\newtheorem{remark}[theorem]{Remark}

\newcommand{\tr}{\mathbb{R}}

\newcommand{\sn}{ S^{n-1}}

\title[A very short unified proof of Dar and log-BM]
{\Large A very short unified proof of  Dar's conjecture   and 
Log-Brunn--Minkowski   in the plane }

\author{\large     Dongmeng Xi}
\address{Department of Mathematics,
Shanghai University,
Shanghai 200444, China}
\email{xi\_dongmeng@shu.edu.cn; dongmeng.xi.math@gmail.com}
 

\begin{document}

\begin{abstract}
This short note gives a very short unified proof of both Dar's conjecture and the log-Brunn--Minkowski inequality in the plane. A planar $L_p$-Brunn-Minkowski inequality for $p\in [-\infty,0]$ is established and gives the two key  inequalities at $p=0$ and $p=-\infty$. Motivated by this new proof, the author   proposes a related  conjecture in $\tr^3$,  to attack the two conjectures in higher dimensions.
\end{abstract}

\maketitle

\section{Introduction}

Dar's conjecture \cite{Dar} states that, for convex bodies
$K,L\subset\tr^n$,
\begin{equation}\label{d1}
 |K+L|^{\frac1n}
 \geq M(K,L)^{\frac1n}
 +\frac{|K|^{\frac1n}|L|^{\frac1n}}{M(K,L)^{\frac1n}},
\end{equation}
where $|\cdot|$ denotes volume, $K+L$ is the Minkowski sum, and
\[
 M(K,L)=\max_{x\in\tr^n}|K\cap(x+L)|.
\]

The log-Brunn--Minkowski conjecture (Log-BM) \cite{BLYZ} states that, for
origin-symmetric convex bodies $K,L\subset\tr^n$ and
$0\leq\lambda\leq1$,
\begin{equation}\label{Lg}
 \left|(1-\lambda)\cdot K+_0\lambda\cdot L\right|
 \geq |K|^{1-\lambda}|L|^\lambda,
\end{equation}
where 
\[
 (1-\lambda)\cdot K+_0\lambda\cdot L
 = 
 \left\{x\in\tr^n:x\cdot v
 \leq h_K(v)^{1-\lambda}h_L(v)^\lambda, \quad \forall v\in S^{n-1} \right\},
\]
and $h_K(v)=\max\{x\cdot v:x\in K\}$ is the support function of $K$.
The log-Brunn--Minkowski conjecture has generated an extensive
literature; see, for example, \cite{BoroczkySurvey} and the references 
therein.  

The \emph{relative inradius} and \emph{relative outradius} of $K$ with
respect to $L$ are defined by
\[
 r(K,L)=\max\{t>0:x+tL\subset K\text{ for some }x\in\tr^n\}
\]
and
\[
 R(K,L)=\min\{t>0:K\subset x+tL\text{ for some }x\in\tr^n\}.
\]
In particular, $r(K,L)=1/R(L,K)$. Throughout this note, whenever there
is no confusion, we simply write $r=r(K,L)$ and $R=R(K,L)$. Two convex
bodies $K$ and $L$ are said to be at a \emph{dilation position} if
\begin{equation}\label{dc}
 o\in K\cap L,
 \qquad rL\subset K\subset RL.
\end{equation}
Every pair can be put at a dilation position after suitable translations
\cite[Lemma~2.1]{XL}. Two origin-symmetric convex bodies are already at such a position.

In dimension 2, the mixed area $V(K,L)$ is
defined by
\[
 |K+tL|=|K|+2tV(K,L)+t^2|L|.
\]
The relative Bonnesen inequality is
\begin{equation}\label{bo1}
 V(K,L) \geq \frac1{t}|K|+ t|L|,
 \qquad r\leq t\leq R.
\end{equation}
Together with the above formula, it gives
\begin{equation}\label{bo2}
 |K+L|\geq\left(1+\frac1t\right)|K|+(1+t)|L|,
 \qquad r\leq t\leq R.
\end{equation}

Xi and Leng \cite{XL} proved both \eqref{d1} and \eqref{Lg} for planar convex bodies that are at a dilation position. The crucial part of their proof was the following alternative inequality:
\[
 \text{either}\quad R|K\cap L|\geq|K|,
 \qquad\text{or}\quad |K\cap L|\geq r|L|,
\]
or, equivalently,
\begin{equation}\label{kcl}
 |K\cap L|\geq
 \min\left\{\frac{|K|}{R},r|L|\right\}.
\end{equation}
Indeed, either $|K|/|K\cap L|$ or $|K\cap L|/|L|$ then belongs to
$[r,R]$. Substituting the corresponding value of $t$ into \eqref{bo2}
gives
\be\label{dar2d}
 |K+L|^{\frac12}
 \geq |K\cap L|^{\frac12}
 +\frac{|K|^{\frac12}|L|^{\frac12}}{|K\cap L|^{\frac12}}.
\ee
Since $|K\cap L|\leq M(K,L)\leq\sqrt{|K||L|}$ and
$s^{1/2}+(|K||L|/s)^{1/2}$ is decreasing for
$0<s\leq\sqrt{|K||L|}$, this implies \eqref{d1}. Thus,   Dar's conjecture in the plane was reduced to
to \eqref{kcl}. The original proof of \eqref{kcl} used a detailed
analysis of the boundary arcs of $K$ and $L$. The same paper also proved  the following general Log-BM for two bodies at a dilation
position:
\be \label{lbm2d}  \left|\frac12\cdot K+_0\frac12\cdot L\right| \ge |K|^{\frac12} \cdot |L|^{\frac12}
\ee 

The purpose of this note is to give a very short proof which naturally
unifies the two inequalities and makes their close relation explicit. Clearly, it suffices to consider the case $\lambda = 1/2$ when dealing with the Log-BM.

For $p\in [-\infty,0]$, the  {\it $L_p$-mean} is defined by
\[
 \frac12\cdot K+_p\frac12\cdot L
 = 
 \left\{x\in\tr^n:x\cdot v\leq
 \left(\frac{h_K(v)^p+h_L(v)^p}{2}\right)^{1/p}, \quad \forall v\in \sn \right\};
\]
at $p=0$, it is understood as the log-combination $ \frac12\cdot K+_0\frac12\cdot L$, and at $p=-\infty$
we put
\[
 \frac12\cdot K+_{-\infty}\frac12\cdot L=K\cap L.
\]

Our main result is the following planar $L_p$-Brunn-Minkowski-type inequality for non-positive $p$, which unifies the two conjectures in the plane. 
By continuity, for $t>0$, we denote  
\[
 \left(\frac{1+t^p}{2}\right)^{2/p}
 =\begin{cases}
 t,&p=0,\\
 \min\{t,1\}^2,&p=-\infty.
 \end{cases}
\]

\begin{theorem}\label{main}
Let $K,L\subset\tr^2$ be at a dilation position. For every
$p\in[-\infty,0]$ and for every $\mu>0$,
\begin{equation}\label{Lp}
 \left|\frac12\cdot K+_p\frac12\cdot L\right|
 \geq\min_{t\in\{r,R\}}
 \left\{\left(\frac{1+t^p}{2}\right)^{2/p}
 \frac{|K|/t+\mu|L|}{t+\mu}\right\}.
\end{equation}
Equality holds if $K$ and $L$ are dilates; when $\mu = 1$ and $|K|=|L|$, equality holds if $K$ and $L$ are dilates or parallelograms with parallel sides.

In particular, 
\begin{itemize}
\item	when  $p=0$ and taking
$\mu=\sqrt{|K|/|L|}$ we obtain the Log-BM \eqref{lbm2d};
\item when $p=-\infty$, $r<1<R$, and taking $ \mu= {r(R-1)}/{(1-r)} $, we obtain \eqref{kcl}. 
\end{itemize}
\end{theorem}
 
As explained above,   Dar's
conjecture follows from \eqref{kcl}. When $R\le 1$ or $r\ge 1$,  \eqref{kcl} follows trivially. Thus, Dar's conjecture is contained in the $-\infty$ case, and hence Theorem \ref{main} unifies both Dar's conjecture and Log-BM in the plane.

The inequality \eqref{Lp} is proved in three steps in Section 2. Motivated by the essential idea in the unified proof, in Section 3,  the author  also introduces the notion of generalized radii and proposes a related formal conjecture in $\tr^3$, in order to attack both the two conjectures in higher dimensions.  

\section{A very short unified proof}

Since this is a short note, we do not recall the general definition of
mixed volumes. In the plane, we use  
\[
 V(K,L)=\frac12\int_{S^1}h_L(v)\,dS_K(v) = \frac12\int_{S^1}h_K(v)\,dS_L(v),\] 
where $S_K$ is the surface area measure of $K$, and the Aleksandrov-Fenchel inequality
\[
 V(K,L)^2\geq |K||L|.
\]
All basic concepts and notation concerning convex bodies and mixed volumes can be found in \cite{Schneider}.

Let $K$ and $L$ be as in Theorem~\ref{main}. For $p\le 0$, write
\[
 M_p=\frac12\cdot K+_p\frac12\cdot L.
\]
We only do the case that the origin is in the interiors of $K$ and $L$;
the boundary case follows by approximation.

\medskip
\noindent\textit{Step 1.}
For $\mu>0$, the function
\[
 s\longmapsto
 \frac{s+\mu}{((1+s^p)/2)^{1/p}}
\]
first decreases and then increases on $(0,\infty)$, since the sign of its
logarithmic derivative is the sign of $1-\mu s^{p-1}$. Since
$r\leq h_K/h_L\leq R$, for $S_{M_p}$-almost every $v\in S^1$,
\[
 \frac{h_{K+\mu L}(v)}{h_{M_p}(v)}
 \leq\max_{t\in\{r,R\}}
 \frac{t+\mu}{((1+t^p)/2)^{1/p}}.
\]
Choose $t\in\{r,R\}$ at which this maximum is attained.

\medskip
\noindent\textit{Step 2.}
The integral formula for mixed areas and the mixed-area inequality give
\[
 \frac{t+\mu}{((1+t^p)/2)^{1/p}}|M_p|
 \geq V(M_p,K+\mu L)
 \geq |M_p|^{1/2}|K+\mu L|^{1/2}.
\]

\medskip
\noindent\textit{Step 3.}
The   Bonnesen inequality \eqref{bo1} gives
\[
\begin{aligned}
 |K+\mu L|
 &=|K|+2\mu V(K,L)+\mu^2|L|\\
 &\geq(t+\mu)\left(\frac{|K|}{t}+\mu|L|\right).
\end{aligned}
\]
Combining the last two inequalities, we obtain
\[
 |M_p|\geq
 \left(\frac{1+t^p}{2}\right)^{2/p}
 \frac{|K|/t+\mu|L|}{t+\mu}.
\]
Since $t$ is either $r$ or $R$, this proves \eqref{Lp} for $p<0$.
Letting $p\to0^-$ and $p\to-\infty$ proves the two endpoint cases.

For completeness, we write down the two short calculations. At $p=0$, for
every $t>0$ and $\mu=\sqrt{|K|/|L|}$,
\[
 t\frac{|K|/t+\mu|L|}{t+\mu}=\sqrt{|K||L|}.
\]
Thus \eqref{lbm2d} follows. 
Note  this inequality holds for every pair of  planar 
convex bodies that are at a dilation position, and it is proved in \cite{XL} that 
  $(1-\lambda)\cdot K+_0\lambda \cdot L$ and $K$ are also at a dilation position for all $\lambda\in[0,1]$.  Then  continuity gives \eqref{Lg} for every $0\leq\lambda\leq1$. 

For $p=-\infty$, the two terms on the right-hand side of \eqref{Lp}, with
$\mu=r(R-1)/(1-r)$, are
\[
 \frac{1-r}{R-r}\frac{|K|}{R}
 +\frac{R-1}{R-r}r|L|
\]
and
\[
 \frac{R(1-r)}{R-r}\frac{|K|}{R}
 +\frac{r(R-1)}{R-r}r|L|.
\]
Both are convex combinations of $|K|/R$ and $r|L|$, and hence give
\eqref{kcl}.

The complete equality characterization will not be discussed here. However, for sufficient condition described in the main theorem, it is merely simple computation.

\begin{remark}
It may be a little difficult to see the main idea of the unified proof of Theorem \ref{main}. The author suggests first assuming $|K|=|L|=\mu=1$ and $p=0$, and then following the three-step proof.
\end{remark}

\section{Generalized   radii and a related conjecture}

The unified planar proof   suggests that, in order to study the two conjectures in higher dimensions, one should look for several relative radii which strengthen the Aleksandrov--Fenchel inequalities in the same way as $r$ and $R$ do in
the plane. 
We first explain the rough idea in $\tr^3$.

We first look for
$t_1,t_2$ such that
\begin{equation}\label{r3}
\begin{aligned}
 3V(K,L,L)&\geq
 \frac{|K|}{t_1t_2}+(t_1+t_2)|L|, \qquad 
 3V(K,K,L)&\geq
 \left(\frac1{t_1}+\frac1{t_2}\right)|K|+t_1t_2|L|.
\end{aligned}
\end{equation}

\vskip 7pt

\noindent{\bf Generalized radii.} 
Define the {\it family of pairs of  generalized radii} by
\be \label{3radii}  \mathcal B_3(K,L):= \left\{ (t_1,t_2)\in(0,\infty)^2 \ \middle|\  \quad {\rm both~inequalities~in~}\eqref{r3} ~ {\rm hold} \right\}. \ee 

\vskip 7pt

In $\tr^2$,   $\mathcal B_2(K,L)$
can be defined as the family of $t$ such that \eqref{bo1} holds, and we see that $[r,R]$ is contained in $\mathcal B_2(K,L)$.  The definition also gives an analogue of \eqref{bo2}:
\[
 |K+L|\geq (1+t_1)(1+t_2)\left ( \frac{|K|}{t_1t_2  }+ |L|\right ),
\]
whenever $(t_1,t_2)\in\mathcal B_3(K,L)$. 

Note that the Aleksandrov--Fenchel inequalities imply 
\[ \frac{|K|}{V(K,K,L)}\cdot |L| \le V(K,L,L), \quad {\rm and} \quad  \frac{V(K,K,L)}{V(K,L,L)}\cdot |L| \le V(K,L,L),\]
and hence 
\[
 \left(\frac{|K|}{V(K,K,L)},
 \frac{V(K,K,L)}{V(K,L,L)} \right)\in\mathcal B_3(K,L).
\]
This implies that   $\mathcal B_3(K,L)$ is nonempty.   

Now, we propose our conjecture. 
Let $K,L\subset\mathbb R^3$ be origin-symmetric convex bodies with $|K|=|L|$, and set as usual: 
\[
M_0=\frac12\cdot K+_0\frac12\cdot L,
\qquad
M_1=\frac{K+L}{2}.
\]
The planar results suggest that, when passing from
the original pair $(K,L)$ to $(M_0,M_1)$, the generalized
radii $(t_1,t_2)$ should be replaced by
\be \label{m0m1radii}
\lambda_1 = \frac{1}{2}\left(\sqrt{t_1}+\frac{1}{\sqrt{t_1}}\right), \quad \lambda_2 = \frac{1}{2}\left(\sqrt{t_2}+\frac{1}{\sqrt{t_2}}\right).
\ee   

\vskip 5pt 

\noindent{\bf Conjecture.}
There is  a  pair of generalized radii
\[
(t_1,t_2)\in\mathcal B_3(K,L)
\]
such that both  inequalities \eqref{conj1} and \eqref{conj2} hold:
\be  \label{conj1} 
3V(M_0,M_0,M_1)
\ge
(\lambda_1+\lambda_2)|M_0|
+\frac{|M_1|}{\lambda_1\lambda_2}
\ee 
and
\be \label{conj2}
3V(M_0,M_0,M_1)
\le
(\lambda_1+\lambda_2)|M_0|
+\frac{|M_0||M_1|}{\lambda_1\lambda_2\sqrt{|K||L|}}.
\ee  
Here $\lambda_i$ is given by \eqref{m0m1radii}. 

\vskip 5pt 

If such a pair can be found, then comparing the  two inequalities gives \(|M_0|\ge \sqrt{|K||L|}\).
 
We briefly explain the  motivation for this conjecture. The conjectured inequality \eqref{conj1} has the same form as
the first inequality in \eqref{r3} and is intended to replace
the Aleksandrov--Fenchel inequality used in the unified planar proof.  The other conjectured inequality \eqref{conj2} is also motivated by the two-dimensional unified proof: In dimension two, if \(|K|=|L|\), take
\[
t=\max\{R,1/r\},
\qquad
\lambda=\frac12\left(\sqrt t+\frac1{\sqrt t}\right).
\]
What we proved in Section 2 implies
\[
2V(M_0,M_1)
\le2\lambda|M_0|
\le\lambda|M_0|
+\frac{|M_0||M_1|}{\lambda\sqrt{|K||L|}},
\]
where Bonnesen's inequality \eqref{bo2} is used once again to derive 
\[ |M_1| \ge \lambda^2 \sqrt{|K||L|}.\]
  
\vskip 20pt

\noindent\textbf{Disclosure.}
The    proof, the notion of generalized radii, and  the proposal of the conjecture  in Section 3 are completely human-made and  due to the author. The proof has been shared with several mathematicians over the past three years.

\end{document}